\documentclass[a4paper,reqno,12pt,oneside]{amsart}
\usepackage{amsmath,amsrefs,geometry}
 \usepackage{color}

\numberwithin{equation}{section}

\renewcommand{\(}{\left(}
\renewcommand{\)}{\right)}

\newtheorem{theorem}{Theorem}[section]

\newtheorem{remark}[theorem]{Remark}

\begin{document}

\date{\today}

\title[Some remarks concerning symmetry-breaking for the Ginzburg-Landau equation]{Some remarks concerning symmetry-breaking for the Ginzburg-Landau equation}

\author{Pierpaolo Esposito}
\address{Pierpaolo Esposito, Dipartimento di Matematica, Universit\`a degli Studi ``Roma Tre'',\,Largo S. Leonardo Murialdo 1, 00146 Roma, Italy}
\email{esposito@mat.uniroma3.it}

\begin{abstract} The correlation term, introduced in \cite{OvSi3} to describe the interaction between very far apart vortices, governs symmetry-breaking for the Ginzburg-Landau equation in $\mathbb{R}^2$ or bounded domains. It is a homogeneous function of degree $(-2)$, and then for $\frac{2\pi}{N}-$symmetric vortex configurations can be expressed in terms of the so-called correlation coefficient. Ovchinnikov and Sigal \cite{OvSi3} have computed it in few cases and conjectured its value to be an integer multiple of $\frac{\pi}{4}$. We will disprove this conjecture by showing that the correlation coefficient always vanishes, and will discuss some of its consequences.
\end{abstract}

\maketitle

\medskip


\section{Introduction}
The Ginzburg-Landau theory is a very popular model in super-conductivity \cite{GiLa}. Stationary states are described by complex-valued solutions $u$ of the planar equation
$$-\Delta u  =k^2 u(1-|u|^2),$$
where $k>0$ is the Ginzburg-Landau parameter. The condensate wave function $u$ describes the superconductive regime in the sample by simply interpreting $|u|^2$ as the density of Cooper electrons pairs. The zeroes of $u$, where the normal state is restored, are called vortices. The parameter $k$ depends on the physical properties of the material and distinguishes between Type I superconductors $k<\frac{1}{\sqrt 2}$ (in this normalization of constants) and Type II superconductors $k>\frac{1}{\sqrt 2}$. 

\medskip \noindent In the entire plane $\mathbb{R}^2$ the parameter $k$ does not play any role, as we can reduce to the case $k=1$ by simply changing $u$ into $u(\frac{x}{k})$. Supplemented by the correct asymptotic behavior at infinity, the Ginzburg-Landau equation now reads as
\begin{equation}
\left\{\begin{array}{ll}
-\Delta U  =U(1-|U|^2)& \hbox{in  }\mathbb{R}^2\\
|U|\to 1 \hbox{ as }|x|\to \infty. &
\end{array} \right. \label{Q}
\end{equation}
The condition $|U|\to 1$ as $|x|\to \infty$ allows to define the (topological) degree $\hbox{deg }U$ of $U$ as the winding number of $U$ at $\infty$: 
$$\hbox{deg }U=\frac{1}{2\pi} \int_{|x|=R} d(arg\, U),$$
where $R>0$ is chosen large so that $|U|\geq \frac{1}{2}$ in $\mathbb{R}^2\setminus B_R(0)$.
Given $n \in \mathbb{Z}$, the only known solution of (\ref{Q}) with $\hbox{deg }U=n$ is the ``radially symmetric" one $U_n(x)=S_n(|x|) (\frac{x}{|x|})^n$ (in complex notations with $x \in \mathbb{C}$), where $S_n$ is the solution of the following ODE:
$$\left\{ \begin{array}{ll} \ddot S_n+\frac{1}{r}\dot S_n-\frac{n^2}{r^2}S_n+S_n(1-S_n^2)=0&\hbox{in }(0,+\infty)\\
S_n(0)=0\:,\:\:\displaystyle \lim_{r \to+\infty}S_n=1. & \end{array}\right. $$
Existence and uniqueness of $S_n$ is shown in \cites{HeHe,CEQ}. Moreover, the solution $U_n$ is stable for $|n|\leq 1$ and unstable for $|n|>1$ \cite{OvSi1}. When $n=\pm 1$, the solution $U_{\pm1}$ is unique, modulo translations and rotations, in the class of functions $U$ with $\hbox{deg }U=\pm 1$ and $\int_{\mathbb{R}^2}(|U|^2-1)^2 dx<+\infty$ \cite{Mir}. 

\medskip \noindent One of the open problems (Problem 1)-- that Brezis-Merle-Rivi\`ere raise out in \cite{BMR}-- concerns the existence of solutions $U$ of (\ref{Q}) with $\hbox{deg }U=n$, $|n|>1$, which are not ``radially symmetric" around any point. So far there is no rigorous answer, but a strategy to find them has been proposed in \cite{OvSi2}. Formally, a solution $U$ of (\ref{Q}) is a critical point of the functional
$$\mathcal{E}(\Psi)=\frac{1}{2} \int_{\mathbb{R}^2} |\nabla \Psi|^2dx+\frac{1}{4} \int_{\mathbb{R}^2} (|\Psi|^2-1)^2 dx.$$
Since $\mathcal{E}(\Psi)=+\infty$ for any $C^1-$map $\Psi$ so that $|\Psi| \to 1$ as $|x|\to +\infty$ and $\hbox{deg }(\Psi)\not=0$, Ovchinnikov and Sigal \cite{OvSi1} have proposed to correct $\mathcal{E}$ into
$$\mathcal{E}_{\hbox{ren}}(\Psi)=\int_{\mathbb{R}^2} \(\frac{1}{2} |\nabla \Psi|^2-\frac{(\hbox{deg }\Psi)^2}{|x|^2}\chi+\frac{1}{4}(|\Psi|^2-1)^2\) dx,$$
 where $\chi$ is a smooth cut-off function with $\chi=0$ when $|x|\leq R$ and $\chi=1$ when $|x|\geq R+R^{-1}$, and $R>>1$ is given. Given a vortex configuration $(\underline a,\underline n)=(a_1,\dots,a_K,n_1,\dots,n_K)$, a $C^1-$map $\Psi$ so that $|\Psi| \to 1$ as $|x|\to +\infty$ has vortex configuration $(\underline a,\underline n)$ if $a_1,\dots,a_K$ are the only zeroes of $\Psi$ with local indices $n_1,\dots,n_K$, denoted for short as $\hbox{conf }\Psi=(\underline a,\underline n)$. Given $\underline n_0$, Ovchinnikov and Sigal \cite{OvSi2} introduce the ``intervortex energy" $E$ given by 
 $$E(\underline a)=\inf \{ \mathcal{E}_{\hbox{ren}} (\Psi):\, \hbox{conf }\Psi=(\underline a,\underline n_0) \},$$
and conjecture that $\underline a_0$ is a critical point of $E$ if and only if there is a minimizer $U$ for $E(\underline a_0)$, yielding to a solution of (\ref{Q}) with $\hbox{conf }U=(\underline a_0,\underline n_0)$ which is not ``radially symmetric" around any point by construction. Letting $d_{\underline a}=\min\limits_{i\not= j}|a_i-a_j|$, the following asymptotic expression is established \cite{OvSi2}:
\begin{equation} \label{asymp1}
E(\underline a)=\sum_{j=1}^K \mathcal{E}_{\hbox{ren}} (U_{n_i})+H(\frac{\underline a}{R})+\hbox{Rem} 
\end{equation}
with $\hbox{Rem} =O(d_{\underline a}^{-1})$ as $d_{\underline a}\to +\infty$, where $H(\underline a)=-\pi \sum \limits_{i \not= j}n_i n_j \ln |a_i-a_j|$ is the energy of the vortex pairs interactions. When $\nabla H(\underline a)=0$, the estimate in (\ref{asymp1}) improves up to $\hbox{Rem} =O(d_{\underline a}^{-2})$.

\medskip \noindent When $\nabla H(\underline a)=0$ (a so-called forceless vortex configuration), by choosing refined test functions the asymptotic expression \eqref{asymp1} is improved \cite{OvSi3} in the form of the following upper bound:
\begin{equation} \label{asymp2}
E(\underline a) \leq \sum_{j=1}^K \mathcal{E}_{\hbox{ren}} (U_{n_i})+H(\frac{\underline a}{R})-A(\underline a)+
\hbox{Rem} 
\end{equation}
with $\hbox{Rem}=O(d_{\underline a}^{-2} +R^{-2})$ as $d_{\underline a}\to +\infty$, where the correlation term $A(\underline a)$ is a homogeneous function of degree $(-2)$ given as
$$A(\underline a)=\frac{1}{4} \int_{\mathbb{R}^2}\left[ |\sum_{j=1}^K \nabla \varphi_j |^4-\sum_{j=1}^K |\nabla \varphi_j|^4 \right],$$
with $\varphi_j(x)=n_j \theta(x-a_j)$, $j=1,\dots,K$, and $\theta(x)$ the polar angle of $x \in \mathbb{R}^2$.

\medskip \noindent To push further the analysis, in \cite{OvSi3} the attention is restricted to symmetric vortex configurations in order to reduce the number of independent variables in $E(\underline a)$. In particular, the simplest $\frac{2\pi}{N}-$symmetric vortex configurations $(\underline a, \underline n)$ (which are invariant under $\frac{2\pi}{N}-$rotations and reflections w.r.t. the real axis) have the form: $a_0=0$, $a_1,\dots,a_N$ are the vertices of a regular $N-$polygon with $a_1=1$ and $n_1=\dots=n_N=m$. We impose also the forceless condition $\nabla H(\underline a)=0$, which simply reads as $n_0=-\frac{N-1}{2}m$.\\
Since $|a_1|=\dots=|a_N|$, the only variable is the size $a=|a_1|$ of the polygon, and then the intervortex energy will be in the form $E(a)$. Since $A(\underline a)$ is homogeneous of degree $-2$, we have that $A(\underline a)=\frac{A_0}{a^2}$, where 
\begin{equation} \label{corrcoeff} A_0:=A(1,e^{\frac{2\pi i}{N}},\dots, e^{\frac{2\pi i(N-1)}{N}})
\end{equation}
is the correlation coefficient for given $n_0= -\frac{N-1}{2}m$ and $n_1=\dots=n_N=m$. In \cite{OvSi3} the existence of c.p.'s of $E(a)$ is shown for the cases $(N,m)=(2,2)$ and $(N,m)=(4,2)$ by comparing $E(a)$ for $a$ small and large, and using the positive sign of $A_0$  (the correlation coefficient has value $8 \pi$ and $80 \pi$, respectively). It is also conjectured \cite{OvSi3} that the correlation coefficient has values which are integer multiples of $\frac{\pi}{4}$. With a long but tricky computation, in the next section we will disprove such a conjecture by showing
\begin{theorem} The correlation coefficient in (\ref{corrcoeff}) always vanishes: $A_0=0$, for all $N\geq 2$ and $m \in \mathbb{Z}$.
\end{theorem}

\medskip \noindent Beside the role of $A_0$ in symmetry-breaking phenomena for (\ref{Q}) in $\mathbb{R}^2$, as already discussed, we will also explain its connection with the Ginzburg-Landau equation 
\begin{equation}
\left\{\begin{array}{ll}
-\Delta u  =k^2 u(1-|u|^2)& \hbox{in  }\Omega\\
u=g &\hbox{on }\partial \Omega
\end{array} \right. \label{QQ}
\end{equation}
on a bounded domain $\Omega$ for strongly Type II superconductors $k \to +\infty$, where $g:\partial \Omega \to S^1$ is a smooth map. 

\medskip \noindent The energy functional for \eqref{QQ}
$$E_k(u)=\frac{1}{2}\int_\Omega |\nabla u|^2+\frac{k^2}{4} \int_\Omega (1-|u|^2)^2$$
has always a minimizer $\bar u_k$ in the space $H=\{u\in H^1(\Omega,\mathbb{C}):\:u=g \hbox{ on }\partial \Omega \}$. When $d=\hbox{deg }g \not= 0$, by \cites{BBH,Str1,Str2} we know that on simply connected domains $\bar u_k$ has exactly $|d|$ simple zeroes $a_1,\dots,a_{|d|}$ for $k$ large, where $(a_1,\dots,a_{|d|})$ is a critical point for a suitable ``renormalized energy" $W(a_1,\dots,a_{|d|})$. The symmetry-breaking phenomenon here takes place, driven by an external mechanism like the boundary condition that forces the confinement of vortices in some equilibrium configuration. A similar result does hold \cite{BBH} on star-shaped domains for any solutions sequence $u_k$ of (\ref{QQ}). Near any vortex $a_i$, the function $u(\frac{x}{k}+a_i)$ behaves like $U_{n_i}(x)$.

\medskip \noindent Once the asymptotic behavior is well understood, a natural question concerns the construction of such solutions for any given c.p. $(a_1,\dots,a_K)$ of $W$, and a positive answer has been given by a heat-flow approach \cites{Lin1,LiLi}, by topological methods \cite{AlBe} and by perturbative methods \cites{MMM2,PaRi} in case $n_1=\dots=n_K= \pm 1$. In \cite{PaRi}, page $12$, it is presented as an open problem to know whether or not there are solutions having
vortices collapsing as $k \to \infty$, the simplest situation being problem (\ref{QQ}) on the unit
ball $B$ with boundary value $g_0=\frac{x^2}{|x|^2}$:
\begin{equation}
\left\{\begin{array}{ll}
-\Delta u  =k^2 u(1-|u|^2)& \hbox{in  }B\\
u=g_0 &\hbox{on }\partial B.
\end{array} \right. \label{P}
\end{equation}
It is conjectured the existence of solutions to (\ref{P}) having a
vortex of degree $-1$ at the origin $a_0=0$ and three vortices of degrees
$+1$ at the vertices $l a_j$, $a_j=e^{\frac{2\pi i}{3} (j-1)}$ for $j=1,2,3$, of a small ($l<<1$) equilateral triangle centered at $0$. This vortex configuration is $\frac{2\pi}{3}-$symmetric, forceless and has ``renormalized energy" 
\begin{equation} \label{W}
W(l)=-6\pi \ln 3-6\pi \ln (1-l^6)+O(l^9)\:,\quad l>0.
\end{equation}

\medskip \noindent In collaboration with J. Wei, we were working on this problem. Inspired by \cite{MMM2}, we were aiming to use a reduction argument of Lyapunov-Schmidt type, starting from the approximating solutions $U_k$ for (\ref{P}) given by
$$U_k(x)=e^{i \varphi_k(x)}U_{-1}(k x) \prod_{j=1}^3 U_{1}\(k(x-l e^{\frac{2\pi i}{3} (j-1)})\)$$
with $l \to 0$ and $lk \to +\infty$, where the function $\varphi_k$ is an harmonic function so that $U_k \mid_{\partial B}=g_0$. The interaction due to the collapsing of three vortices onto $0$ gives at main order a term $(lk)^{-2}$ with the plus sign, i.e. for some $J_0>0$ there holds the energy expansion
\begin{eqnarray} \label{expansionenergy}
E_k(U_k)&=&4\pi \ln k+I+\frac{1}{2} W(l)+J_0 (lk)^{-2}+o((lk)^{-2})\nonumber \\
&=& 4\pi \ln k+I-3 \pi \ln 3+3\pi l^6+J_0 (lk)^{-2}+o\(l^6+(lk)^{-2}\),
\end{eqnarray}
in view of \eqref{W}. The aim is to construct a solution $u_k$ in the form $ U_k [\eta(1+\psi)+(1-\eta)e^{\psi}]$, where $\psi=\psi(k)$ is a remainder term small in a weighted $L^\infty(B)-$norm and $l=l(k)$ as $k \to +\infty$. The function $\eta$ is a smooth cut-off function with $\eta=1$ in $\cup_{j=0}^3 B_{1/k}(la_j)$ and $\eta=0$ in $B \setminus \cup_{j=0}^3 B_{2/k}(la_j)$. The function $\psi=\psi(k)$ is found thanks to the solvability theory (up to a finite-dimensional kernel) of the linearized operator for \eqref{P} at $U_k$ as $l\to 0$ and $l k\to +\infty$, and by the Lyapunov-Schimdt reduction the existence of $l(k)$ follows  as a c.p. of 
$$\tilde E_k:=E_k\( U_k [\eta(1+\psi(k))+(1-\eta)e^{\psi(k)}] \).$$
If $U_k$ is sufficiently good as an approximating solution of (\ref{P}), we have that $\tilde E_k=E_k(U_k)+o((lk)^{-2})$. Since $3\pi l^6+J_0 (lk)^{-2}$ has always a minimum point of order $k^{-\frac{1}{4}}$ as $k \to +\infty$, by \eqref{expansionenergy} we get the existence of $l=l(k)$ in view of the persistence of minimum points under small perturbations.

\medskip \noindent Unfortunately, this is not the case. Pushing further the analysis, we were able to identify the leading term $\psi_0=\psi_0(k)$ of $\psi=\psi(k)$, and compute its contribution into the energy expansion, yielding to a correction in the form:
\begin{equation} \label{expansionenergybis}
\tilde E_k=4\pi \ln k+I+\frac{1}{2} W(l)+J_1 (lk)^{-2}+o((lk)^{-2}).
\end{equation}
By \eqref{W} and \eqref{expansionenergybis} a c.p. $l(k)$ of $\tilde E_k$ always exists provided $J_1>0$. First numerically, and then rigorously, we were disappointed to find that $J_1=0$. 

\medskip \noindent Later on, we realized that $-J_1$ is exactly the correlation coefficient $A_0$ in (\ref{corrcoeff}) (with $N=3$ and $m=1$) introduced by Ovchinnikov and Sigal \cite{OvSi3}. If $u$ is a solution of (\ref{P}) with vortices $a_0=0$ and $l a_j$, $a_j=e^{\frac{2\pi i}{3} (j-1)}$ for $j=1,2,3$, with $n_0=-1$ and $n_1=n_2=n_3=1$, then the function $U(x)=u(\frac{x}{k})$ does solve
\begin{equation}
\left\{\begin{array}{ll}
-\Delta U  =U(1-|U|^2)& \hbox{in  }B_k\\
U=g_0 &\hbox{on }\partial B_k
\end{array} \right. \label{Pbis}
\end{equation}
with vortices $a_0$ and $l k a_j$ of vorticities $n_0=-1$, $n_1=n_2=n_3=1$. Since (\ref{Q}) and (\ref{Pbis}) formally coincide when $k=+\infty$, it is natural to find a correlation term in the energy expansion $\tilde E_k$ in the form $-\frac{A_0}{a^2}=J_1(l k)^{-2}$, where $a=lk$ is the modulus of the $l k a_j$'s for $j=1,2,3$. Even more and not surprisingly, the function $\tilde U_k(\frac{x}{k})$, where
$$U_k [\eta(1+\psi_0(k))+(1-\eta)e^{\psi_0(k)}]$$
is a very good approximating solution for (\ref{P}) which improves the approximation rate of $U_k$, does coincide with the refined test functions used by Ovchinnikov and Sigal \cite{OvSi3} to get the improved upper bound (\ref{asymp2}).

\medskip \noindent In conclusion, the vanishing of the correlation coefficient $A_0$ does not support any conjecture  concerning symmetry-breaking phenomena for (\ref{Q}) or the existence of collapsing vortices for (\ref{P}) when $k \to +\infty$. Higher-order expansions would be needed in their study.

\section{The correlation coefficient} 
\noindent Let $N\geq 2$. Let $a_j=e^{\frac{2\pi i (j-1)}{N}}$, $j=1,\dots,N$, be the $N-$roots of unity, and set $n_j=m \in \mathbb{Z}$ for all $j=1,\dots,N$, $a_0=0$ and $n_0=-\frac{N-1}{2}m$. We aim to compute the correlation coefficient $A_0=A_0(m)$ given in \eqref{corrcoeff}. Since (in complex notation) $\nabla \theta(x)=|x|^{-2}(-x_2,x_1)$ has the same modulus as $\frac{1}{x}=\frac{\bar x}{|x|^2}$, the correlation coefficient takes the form
\begin{equation} \label{correlationcoeffequiv}
A_0=\frac{1}{4} \int_{\mathbb{R}^2}\left[ |\sum_{j=0}^N \frac{n_j}{x-a_j} |^4-\sum_{j=0}^N |\frac{n_j}{x-a_j}|^4 \right].
\end{equation}
Since the integer $m$ comes out as $m^4$ from the expression (\ref{correlationcoeffequiv}), we have that $A_0(m)=m^4 A_0(1)$. Hereafter, we will assume $m=1$ and simply denote $A_0(1)$ as $A_0$.

\medskip \noindent Let us first notice that $A_0$ is not well-defined without further specifications, because the integral function in (\ref{correlationcoeffequiv}) is not integrable near the points $a_j$, $j=0,\dots,N$. Recall that the $N-$roots of unity $a_1,\dots,a_N$ do satisfy the following symmetry properties:
\begin{equation} \label{symmetry}
\sum_{j=1}^N a_j^l=0 \qquad \forall\, |l| \leq N,\, l\not=0,
\end{equation}
as it can be easily deduced by the relation $x^N-1=\prod\limits_{j=1}^N (x-a_j)$. A first application of (\ref{symmetry}) is the validity of
\begin{equation} \label{sum}
\sum_{j=1}^N \frac{1}{x-a_j}=\sum_{j=1}^N \frac{x^{N-1}+a_j x^{N-2}+\dots+a_j^{N-1}}{x^N-1}=\frac{N x^{N-1}}{x^N-1},
\end{equation}
which implies that the integral function in \eqref{correlationcoeffequiv} near $0$ has the form
\begin{eqnarray} \label{near0}
|\sum_{j=0}^N \frac{n_j}{x-a_j} |^4-\sum_{j=0}^N |\frac{n_j}{x-a_j}|^4=- \frac{N(N-1)^3}{2} \hbox{Re}(\frac{x^{N}}{(x^N-1)|x|^4}) +O(1)
\end{eqnarray}
and is not integrable at $0$ when $N=2$. Similarly, setting $\alpha_k(x)=-\frac{N-1}{2x}+\sum\limits_{j=1 \atop j \not= k}^N \frac{1}{x-a_j}$ for $k=1,\dots,N$, near $a_k$ we have that
\begin{eqnarray} \label{nearajpreliminary}
&& |\sum_{j=0}^N \frac{n_j}{x-a_j} |^4-\sum_{j=0}^N |\frac{n_j}{x-a_j}|^4=
\frac{4}{|x-a_k|^4} \hbox{Re}[(x-a_k)\alpha_k(x)]+\frac{2}{|x-a_k|^2}|\alpha_k(x)|^2 \\
&&+\(2 \hbox{Re}\frac{(x-a_k)\alpha_k(x)}{|x-a_k|^2}+|\alpha_k(x)|^2 \)^2-\frac{(N-1)^4}{16|x|^4}-\sum_{j=1 \atop j \not= k}^N \frac{1}{|x-a_j|^4}.
\nonumber
\end{eqnarray}
The function $\alpha_k$ can not be computed explicitly, but we know that
\begin{eqnarray} \label{alphak}
\alpha_k(a_k)&=&-\frac{N-1}{2a_k}+\sum_{j=1 \atop j\not= k}^N \frac{1}{a_k-a_j}=
a_k^{N-1} \(-\frac{N-1}{2}+\sum_{j=2}^N \frac{1}{1-a_j }\)\\
&=&a_k^{N-1} \(-\frac{N-1}{2}+\sum_{j=2}^N \frac{1-\cos \frac{2\pi (j-1)}{N}+i\sin \frac{2\pi (j-1)}{N}}{2(1-\cos \frac{2\pi (j-1)}{N})} \) \nonumber \\
&=&  i a_k^{N-1}  \sum_{j=2}^N \frac{\sin \frac{2\pi (j-1)}{N}}{2(1-\cos \frac{2\pi (j-1)}{N})}=0 \nonumber 
\end{eqnarray}
in view of $\{a_j a_k^{N-1}:\, j=1,\dots,N,\, j\not= k\}=\{a_2,\dots,a_N\}$ and the symmetry of $\{a_1,\dots,a_N\}$  under reflections w.r.t. the real axis. By inserting \eqref{alphak} into \eqref{nearajpreliminary} we deduce that the integral function in \eqref{correlationcoeffequiv} near $a_k$ has the form
\begin{eqnarray} \label{nearaj}
|\sum_{j=0}^N \frac{n_j}{x-a_j} |^4-\sum_{j=0}^N |\frac{n_j}{x-a_j}|^4=
\frac{4}{|x-a_k|^4} \hbox{Re}[\alpha_k'(a_k)(x-a_k)^2]+O(\frac{1}{|x-a_k|})
\end{eqnarray}
and is not integrable at $a_k$ when $\alpha_k'(a_k) \not=0$. Since the (possible) singular term in \eqref{near0}, \eqref{nearaj} has vanishing integrals on circles, the meaning of $A_0$ is in terms of a principal value:
\begin{equation} \label{correlationcoeffequivprecise}
A_0=\frac{1}{4} \lim\limits_{\epsilon \to 0}  \int_{\mathbb{R}^2 \setminus \cup_{k=0}^N B_\epsilon(a_k)}\left[ |\sum_{j=0}^N \frac{n_j}{x-a_j} |^4-\sum_{j=0}^N |\frac{n_j}{x-a_j}|^4 \right].
\end{equation}

\medskip \noindent We would like to compute $A_0$ in polar coordinates, even tough the set $\mathbb{R}^2 \setminus \cup_{k=0}^N B_\epsilon(a_k)$ is not radially symmetric. The key idea is to make the integral function in \eqref{correlationcoeffequivprecise} integrable near any $a_j$, $j=1,\dots,N$, by adding suitable singular terms, in such a way that the integral in \eqref{correlationcoeffequivprecise} will have to be computed just on the radially symmetric set $\mathbb{R}^2 \setminus B_\epsilon(a_0)$. To this aim, it is crucial to compute $\alpha_k'(a_k)$. Arguing as before, we get that
\begin{eqnarray} \label{derivalphak}
\alpha_k'(a_k)&=&\frac{N-1}{2 a_k^2}-\sum_{j=1 \atop j\not= k}^N \frac{1}{(a_k-a_j)^2}=
a_k^{N-2} \(\frac{N-1}{2}-\sum_{j=2}^N \frac{1}{(1-a_j)^2 }\) \nonumber \\
&=&a_k^{N-2} \(\frac{N-1}{2}-\sum_{j=2}^N \frac{(1-\cos \frac{2\pi (j-1)}{N})^2-\sin^2 \frac{2\pi (j-1)}{N}}{4(1-\cos \frac{2\pi (j-1)}{N})^2} \) \nonumber \\
&=&a_k^{N-2} \sum_{j=2}^N \frac{1}{2(1-\cos \frac{2\pi (j-1)}{N})}=a_k^{N-2} \sum_{j=2}^N \frac{1}{|1-a_j|^2}.
\end{eqnarray}
Since there holds $\sum\limits_{j=1}^{N-1} a_k^j=\sum\limits_{j=2}^N a_j=-1$ for all $k=2,\dots,N$ in view of \eqref{symmetry}, we have that
$$\prod\limits_{j=2}^N (z-a_j)=\frac{z^N-1}{z-1}=\sum\limits_{p=0}^{N-1} z^p \,,\quad
\prod\limits_{j=2 \atop j\not=k}^N (z-a_j)=\frac{\sum\limits_{p=0}^{N-1}z^p}{z-a_k}=\sum\limits_{p=0}^{N-2} z^p \sum_{l=0}^{N-2-p} a_k^l,$$
and then
\begin{eqnarray} \label{pkl}
\prod\limits_{j=2}^N (1-a_j)=N \,,\quad
\prod\limits_{j=2 \atop j\not=k}^N (1-a_j)=\sum\limits_{l=0}^{N-2} (N-l-1) a_k^l.
\end{eqnarray}
By \eqref{pkl} we get that
\begin{eqnarray*}
\beta_N &:=& \sum_{j=2}^N \frac{4}{|1-a_j|^2}=\sum_{j=2}^N \frac{4}{N^2} \prod\limits_{k=2 \atop k \not=j }^N |1-a_k|^2 =
\sum_{j=2}^N \frac{4}{N^2} \sum\limits_{l,p=0}^{N-2} (N-l-1)(N-p-1) a_j^{l-p} \\
&=& 4\frac{N-1}{N^2} \sum\limits_{l=1}^{N-1} l^2-
\frac{4}{N^2} \sum\limits_{l,p=1 \atop l\not=p}^{N-1} lp =
\frac{4}{N} \sum\limits_{l=1}^{N-1} l^2-\frac{4}{N^2} (\sum\limits_{l=1}^{N-1} l)^2=
\frac{2(N-1)(2N-1)}{3} -(N-1)^2\\
&=& \frac{N^2-1}{3}
\end{eqnarray*}
in view of \eqref{symmetry}. Since by \eqref{derivalphak} $\alpha_k'(a_k)=\frac{\beta_N}{4} a_k^{N-2}$, by \eqref{nearaj} we have that
$$|\sum_{j=0}^N \frac{n_j}{x-a_j} |^4-\sum_{j=0}^N |\frac{n_j}{x-a_j}|^4-
 \sum_{j=1}^N \hbox{Re}[\frac{\beta_N a_j^2 }{(x-a_j)^2 (1+|x-a_j|^2)} ] \in L^1(\mathbb{R}^2 \setminus \{0\}).$$
Since
$$\lim\limits_{\epsilon \to 0}  \int_{\mathbb{R}^2 \setminus \cup_{k=0}^N B_\epsilon(a_k)} \frac{a_j^2}{(x-a_j)^2 (1+|x-a_j|^2)}=\lim\limits_{\epsilon \to 0}  \int_{\mathbb{R}^2 \setminus B_\epsilon(a_j)} \frac{a_j^2}{(x-a_j)^2 (1+|x-a_j|^2)}=0,$$
we can re-write $A_0$ as
\begin{eqnarray} \label{A}
A_0&=& \frac{1}{4}\lim\limits_{\epsilon \to 0}  \int_{\mathbb{R}^2 \setminus B_\epsilon(0)}\left[ |\frac{(N+1) x^N+(N-1)}{2x(x^N-1)}|^4-\frac{(N-1)^4}{16 |x|^4}-\sum_{j=1}^N \frac{1}{|x-a_j|^4}\right. \nonumber \\
&&-\left. \sum_{j=1}^N \hbox{Re}[\frac{\beta_N a_j^2}{(x-a_j)^2 (1+|x-a_j|^2)} ] \right]\nonumber\\
&=& \frac{1}{4}\lim\limits_{\epsilon \to 0}  \int_{\mathbb{R}^2 \setminus (B_\epsilon(0) \cup  \{1-\epsilon\leq |x|\leq \frac{1}{1-\epsilon}\})} \left[ |\frac{(N+1) x^N+(N-1)}{2x(x^N-1)}|^4-\frac{(N-1)^4}{16 |x|^4}-\sum_{j=1}^N \frac{1}{|x-a_j|^4}\right] \nonumber \\
&&- \frac{1}{4}  \hbox{Re}\left[ \lim\limits_{\epsilon \to 0} \int_{\mathbb{R}^2 \setminus (B_\epsilon(0) \cup  \{1-\epsilon\leq |x|\leq \frac{1}{1-\epsilon}\})} \sum_{j=1}^N \frac{\beta_N a_j^2}{(x-a_j)^2(1+|x-a_j|^2)}\right]=:\frac{1}{4}\hbox{I}-\frac{1}{4}\hbox{II}
\end{eqnarray}
in view of \eqref{sum}.

\medskip \noindent As far as $\hbox{I}$, let us write the following Taylor expansions: for $|x|<1$ there hold
\begin{eqnarray} \label{Taylor1}
\frac{((N+1)x^N+(N-1))^2}{(1-x^N)^2}&=&\((N-1)^2+2(N^2-1)x^N+(N+1)^2 x^{2N} \) \sum_{k \geq 0}(k+1)x^{kN} \nonumber \\
&=&(N-1)^2+\sum_{k\geq 1}4N(kN-1) x^{kN}=\sum_{k\geq 0}c_k x^{kN}
\end{eqnarray}
and
\begin{eqnarray} \label{Taylor2}
\frac{((N-1)x^N+(N+1))^2}{(1-x^N)^2}&=& \((N+1)^2+2(N^2-1) x^N +(N-1)^2 x^{2N}\)\sum_{k\geq 0}(k+1)x^{kN} \nonumber \\
&=& (N+1)^2+\sum_{k\geq 1}4N(kN+1)x^{kN}=\sum_{k\geq 0}d_k x^{kN},
\end{eqnarray}
where $c_k=\max\{4N(kN-1),(N-1)^2 \}$ and $d_k=\max\{4N(kN+1),(N+1)^2\}$. Letting $\epsilon>0$ small, by (\ref{Taylor1})-(\ref{Taylor2}) we have that in polar coordinates (w.r.t. to the origin)  $\hbox{I}$ writes as
\begin{eqnarray*}
\hbox{I}&=& \int_\epsilon^{1-\epsilon} \rho d\rho \int_0^{2\pi}d\theta \left[ \frac{1}{16 \rho^4} |\sum_{k\geq 0} c_k \rho^{kN}e^{ikN \theta}|^2-\frac{(N-1)^4}{16 \rho^4}- \sum_{j=1}^N
|\sum_{k\geq 0}(k+1)a_j^{k(N-1)} \rho^k e^{ik\theta}|^2  \right]\\
&&+\int_{\frac{1}{1-\epsilon}}^\infty \rho d\rho \int_0^{2\pi}d\theta \left[ \frac{1}{16 \rho^4} |\sum_{k\geq 0} d_k \rho^{-k N}e^{-ik N \theta}|^2-\frac{(N-1)^4}{16 \rho^4}-\frac{1}{\rho^4}\sum_{j=1}^N
|\sum_{k\geq 0}(k+1)a_j^k \rho^{-k} e^{-ik\theta}|^2  \right]\\
&&+o_\epsilon(1)
\end{eqnarray*}
with $o_\epsilon(1)\to 0$ as $\epsilon \to 0$, in view of
$$|x-a_j|^{-4}=|a_j^{N-1}x-1|^{-4}=|\sum_{k\geq 0}(k+1)a_j^{k(N-1)}x^k|^2,\quad |1-a_j x|^{-4}=|\sum_{k\geq 0}(k+1)a_j^k x^k|^2$$ for $|x|<1$. By the Parseval's Theorem  we get that
\begin{eqnarray*}
\hbox{I}&=& 2\pi \int_\epsilon^{1-\epsilon} \left[\frac{1}{16} \sum_{k\geq 1} |c_k|^2 \rho^{2kN-3}-N \sum_{k\geq 0}
(k+1)^2 \rho^{2k+1}\right] d \rho\\
&&+2\pi \int_{\frac{1}{1-\epsilon}}^\infty \left[\frac{1}{16} \sum_{k\geq 1} |d_k|^2 \rho^{-2KN-3}+\frac{(N+1)^4-(N-1)^4}{16 \rho^3}-N \sum_{k\geq 0}(k+1)^2 \rho^{-2k-3}\right] d \rho+o_\epsilon(1)\\
&=& 2\pi N \int_0^{1-\epsilon} \left[ N \sum_{k\geq 0} (kN+N-1)^2 \rho^{2kN+2N-3}-\sum_{k\geq 0}
(k+1)^2 \rho^{2k+1}\right] d \rho\\
&&+2 \pi N \int_{\frac{1}{1-\epsilon}}^\infty \left[N \sum_{k\geq 0} (kN+N+1)^2 \rho^{-2kN-2N-3}-\sum_{k\geq 0}(k+1)^2 \rho^{-2k-3}\right] d \rho+N(N^2+1) \frac{\pi}{2}\\
&&+o_\epsilon(1)= 2\pi N \int_0^{1-\epsilon} \left[ N \sum_{k\geq 0} (kN+N-1)^2 \rho^{2kN+2N-3}
+N \sum_{k\geq 0} (kN+N+1)^2 \rho^{2kN+2N+1}\right.\\
&&\left. -2 \sum_{k\geq 0} (k+1)^2 \rho^{2k+1}\right] d \rho+N(N^2+1) \frac{\pi}{2}+o_\epsilon(1)
\end{eqnarray*}
as $\epsilon \to 0$. We compute now the integrals and let $\epsilon \to 0$ to end up with
\begin{eqnarray*}
\hbox{I}&=& 2\pi N \left[\frac{N}{2} \sum_{k\geq 0} (kN+N-1) \rho^{2kN+2N-2}+
\frac{N}{2} \sum_{k\geq 0} (kN+N+1)\rho^{2kN+2N+2} -\sum_{k\geq 0} (k+1) \rho^{2k+2}\right] \Big|_0^1\\
&&+N(N^2+1) \frac{\pi}{2}.
\end{eqnarray*}
Denoting the function inside brackets as $f(\rho)$, we need now to determine the explicit expression of $f(\rho)$ for $\rho<1$:
\begin{eqnarray*}
f(\rho)&=& \frac{N^2}{2}  \rho^{2N-2}(1+\rho^4) \sum_{k\geq 0} (k+1) (\rho^{2N})^k-
\frac{N}{2} \rho^{2N-2}(1-\rho^4) \sum_{k\geq 0} (\rho^{2N})^k 
- \rho^2 \sum_{k\geq 0} (k+1)(\rho^{2})^k \\
&=&\frac{N^2}{2}  \rho^{2N-2}\frac{1+\rho^4}{(1-\rho^{2N})^2}-
\frac{N}{2} \rho^{2N-2}\frac{1-\rho^4}{1-\rho^{2N}} 
- \frac{\rho^2}{(1-\rho^2)^2}\\
&=&\frac{1}{2} \frac{N^2 \rho^{2N-2}(1+\rho^4)-N \rho^{2N-2}(1-\rho^4)(1-\rho^{2N})-2\rho^2 (\sum\limits_{j=0}^{N-1} \rho^{2j})^2}{(1-\rho^{2N})^2},
\end{eqnarray*}
and then by the l'H\^{o}pital's rule we get that
\begin{eqnarray*}
4 N^2 f(1)&=& 2 \lim\limits_{\rho \to 1} \frac{N(N-1) \rho^{N-1}+N (N+1) \rho^{N+1}-2\rho (\sum\limits_{j=0}^{N-1} \rho^j)^2+N \rho^{2N-1}-N \rho^{2N+1}}{(1-\rho)^2}\\
&=&\lim\limits_{\rho \to 1} \frac{-N^2(N-2) \rho^{N-2}-N^2 (N+2) \rho^N+2(\sum\limits_{j=0}^{N-1} \rho^j)^2+4 \rho (\sum\limits_{j=0}^{N-1} \rho^j) (\sum\limits_{j=0}^{N-2} (j+1) \rho^j)}{1-\rho}\\
&&+N \lim\limits_{\rho \to 1} \frac{(2N+1) \rho^{2N}-(2N-1) \rho^{2N-2}-\rho^{N-2}-\rho^N}{1-\rho}=-\frac{N^2 (N^2+5)}{3}.
 \end{eqnarray*}
In conclusion, for $\hbox{I}$ we get the value
\begin{eqnarray} \label{valueI}
\hbox{I}= \frac{\pi}{3} N(N^2-1).
\end{eqnarray}
\begin{remark} In \cite{OvSi3} the value of $A_0$ was computed neglecting the term $\hbox{II}$ in \eqref{A}. By \eqref{valueI} notice that $\frac{m^4}{4}\hbox{I}=\frac{\pi}{12} m^4 N(N^2-1)$ does coincide with $8 \pi$ when $(N,m)=(2,2)$ and $80 \pi$ when $(N,m)=(4,2)$, in agreement with the computations in \cite{OvSi3}.
\end{remark}

\medskip \noindent As far as $\hbox{II}$, let us compute in polar coordinates the value of
\begin{eqnarray*}
\lim\limits_{\epsilon \to 0} \sum_{j=1}^N \int_{\mathbb{R}^2 \setminus (B_\epsilon(0) \cup  \{1-\epsilon\leq |x|\leq \frac{1}{1-\epsilon}\})}  \frac{a_j^2}{(x-a_j)^2(1+|x-a_j|^2)}= \lim\limits_{\epsilon \to 0} \int_{(0,1-\epsilon)\cup (\frac{1}{1-\epsilon},+\infty)} \rho \Gamma(\rho) d\rho,
\end{eqnarray*}
where the function $\Gamma$ is defined in the following way:
\begin{eqnarray*}
\Gamma(\rho)&=& \sum_{j=1}^N \int_0^{2\pi} \frac{a_j^2}{(\rho e^{i\theta}-a_j)^2(2+\rho^2-a_j \rho e^{-i\theta}-a_j^{N-1} \rho e^{i\theta})  }d \theta\\
&=& \frac{i}{\rho} \sum_{j=1}^N a_j^3 \int_\gamma \frac{dw}{(\rho w-a_j)^2(w^2-\frac{2+\rho^2}{\rho}a_j w+a_j^2)},
\end{eqnarray*}
with $\gamma$ the counterclockwise unit circle around the origin. Since
$$w^2-\frac{2+\rho^2}{\rho}a_j w+a_j^2=\left(w-\frac{2+\rho^2}{2\rho}a_j\right)^2+a_j^2\left(1-(\frac{2+\rho^2}{2\rho})^2 \right),$$
observe that $w^2-\frac{2+\rho^2}{\rho}a_j w+a_j^2$ vanishes at $\rho_\pm a_j$, with
$$\rho_\pm=\frac{2+\rho^2}{2\rho} \pm \sqrt{ (\frac{2+\rho^2}{2\rho})^2-1 }$$
satisfying $\rho_-<1<\rho_+$ in view of $\frac{2+\rho^2}{2\rho}\geq \sqrt 2$. Since
$$(\frac{1}{w^2-\frac{2+\rho^2}{\rho} a_j w+a_j^2 })' (\frac{a_j}{\rho})=a_j^{N-3}\rho^5,$$
by the Cauchy's residue Theorem the function $\Gamma(\rho)$ can now be computed explicitly as
\begin{eqnarray*} 
\Gamma(\rho)&=& \frac{i}{\rho^3}\sum_{j=1}^N a_j^3 \int_\gamma \frac{dw}{(w-\frac{a_j}{\rho})^2(w-\rho_-a_j)(w-\rho_+a_j)}\\
&=& 2\pi N \left \{ \begin{array}{ll} (\rho \rho_- -1)^{-2}(\rho \rho_+-\rho \rho_-)^{-1} &\hbox{if }\rho<1\\
(\rho \rho_- -1)^{-2} (\rho \rho_+ - \rho \rho_-)^{-1}-\rho^2 &\hbox{if }\rho>1. \end{array} \right.
\end{eqnarray*}
Since we have that
$$(\rho \rho_- - 1)^2=\frac{1}{4}(\rho^2-\sqrt{\rho^4+4})^2=\frac{1}{2} (\rho^4+2-\rho^2\sqrt{\rho^4+4}),\qquad \rho \rho_+ - \rho \rho_-=\sqrt{\rho^4+4},$$
we get that
$$(\rho \rho_- - 1)^{-2}(\rho \rho_+ - \rho \rho_-)^{-1}=\frac{2}{(\rho^4+2)\sqrt{\rho^4+4}-\rho^2(\rho^4+4)}=\frac{\rho^4+2}{2\sqrt{\rho^4+4}}+\frac{\rho^2}{2},$$
and the expression of $\Gamma(\rho)$ now follows in the form
\begin{eqnarray} \label{Gamma}
\Gamma(\rho)= \pi N \frac{\rho^4+2}{\sqrt{\rho^4+4}}- \pi N \rho^2 +
\left \{ \begin{array}{ll} 2 \pi N \rho^2 &\hbox{if }\rho<1\\
0 &\hbox{if }\rho>1. \end{array} \right.
\end{eqnarray}
Note that
$$\rho \left( \frac{\rho^4+2}{\sqrt{\rho^4+4}}-\rho^2  \right)=\frac{4\rho}{(\rho^4+2) \sqrt{\rho^4+4}+\rho^2(\rho^4+4)}$$
is integrable in $(0,\infty)$, and we have that
\begin{eqnarray} \label{integral1}
&&\int_0^\infty \rho (\frac{\rho^4+2}{\sqrt{\rho^4+4}} -\rho^2)d\rho=\lim_{M \to +\infty} \frac{1}{2} \int_0^M (\frac{s^2+2}{\sqrt{s^2+4}}-s)ds\\
&&=\lim_{M \to +\infty}\left[  \frac{s}{4}\sqrt{s^2+4} \mid_0^M-\frac{M^2}{4} \right]=\lim_{M \to +\infty}\frac{M}{4} (\sqrt{M^2+4}-M)=\frac{1}{2}.\nonumber
\end{eqnarray}
Thanks to (\ref{Gamma})-(\ref{integral1}) we can compute 
\begin{eqnarray*}
\lim\limits_{\epsilon \to 0} \int_{(0,1-\epsilon)\cup (\frac{1}{1-\epsilon},+\infty)} \rho \Gamma(\rho) d\rho&=&
\int_0^{+\infty} \rho \Gamma(\rho) d\rho=\pi N,
\end{eqnarray*}
and for $\hbox{II}$ we get the value
\begin{eqnarray} \label{valueII}
\hbox{II}= \frac{\pi}{3} N (N^2-1) .\end{eqnarray}

\medskip \noindent Finally, inserting (\ref{valueI}) and (\ref{valueII}) into (\ref{A}) we get that the correlation coefficient vanishes: $A_0=0$. Then, there holds $A_0(m)=0$ for all $m \in \mathbb{Z}$, as claimed.


\begin{thebibliography}{99}

\bib{AlBe}{article}{
   author={Almeida, L.},
   author={Bethuel, F.},
   title={Topological methods for the Ginzburg-Landau equations},
   journal={J. Math. Pures Appl. (9)},
   volume={77},
   date={1998},
   number={1},
   pages={1--49},
}

\bib{BBH}{book}{
   author={Bethuel, F.},
   author={Brezis, H.},
   author={H\'elein, F.},
   title={Ginzburg-Landau vortices},
   series={Progress in Nonlinear Differential Equations and their Applications},
   volume={13},
   publisher={Birkh\"auser},
   place={Boston},
   date={1994},
}

\bib{BMR}{article}{
   author={Brezis, H.},
   author={Merle, F.},
   author={Rivi\`ere, T.}
   title={Quantization effects for $-\Delta u=u(1-|u|^2)$ in $\mathbb{R}^2$},
   journal={Arch. Rat. Mech. Anal.},
   volume={126},
   date={1994},
   number={1},
   pages={35--58},
}

\bib{CEQ}{article}{
    author= {Chen, X.},
    author= {Elliott, C.M.},
    author= {Qi, T.},
    title= {Shooting method for vortex solutions of a complex-valued
              Ginzburg-Landau equation},
    journal = {Proc. Roy. Soc. Edinburgh Sect. A},
    volume= {124},
    date= {1994},
    number= {6},
    pages= {1075--1088},
}



\bib{MMM2}{article}{
   author={del Pino, M.},
   author={Kowalczyk, M.},
   author={Musso, M.}
   title={Variational reduction for Ginzburg-Landau vortices},
   journal={J. Funct. Anal.},
   volume={239},
   date={2006},
   number={2},
   pages={497--541},
}

\bib{GiLa}{article}{
   author={Ginzburg, V.L.},
   author={Landau, L.D.},
   title={On the theory of superconductivity},
   journal={J. Exp. Theor. Phys.},
   volume={20},
   date={1950},
   pages={1064--1082},
}

\bib{HeHe}{article}{
   author={Herv\'e, R.-M.},
   author={Herv\'e, M.},
   title={ \'Etude qualitative des solutions r\'eelles d'une \'equation diff\'erentielle li\'ee \`a l'\'equation
de Ginzburg-Landau},
   journal={Ann. Inst. H. Poincar\'e Anal. Non Lin{\'e}aire},
   volume={11},
   date={1994},
   number={4},
   pages={427--440},
}

\bib{Lin1}{article}{
   author={Lin, F.H.},
   title={Solutions of Ginzburg-Landau equations and critical points
              of the renormalized energy},
   journal={Ann. Inst. H. Poincar\'e Anal. Non Lin{\'e}aire},
   volume={12},
   date={1995},
   number={5},
   pages={599--622},
}

\bib{LiLi}{article}{
   author={Lin, F.H.},
   author={Lin, T.-C.},
   title={Minimax solutions of the Ginzburg-Landau equations},
   journal={Selecta Math. (N.S.)},
   volume={3},
   date={1997},
   number={1},
   pages={99--113},
}


\bib{Mir}{article}{
   author={Mironescu, P.},
   title={Local minimizers for the Ginzburg-Landau equation are radially symmetric},
   journal={C. R. Acad. Sci. Paris S\'er. I Math.},
   volume={323},
   date={1996},
   number={6},
   pages={593--598},
}
 
\bib{OvSi1}{book}{
   author={Ovchinnikov, Yu.N.},
   author={Sigal, I.M.},
   title={Ginzburg-Landau equation. I. Static vortices},
   series= {CRM Proc. Lecture Notes, Partial differential equations and their applications
               (Toronto, ON, 1995)},
   pages={199--220},
   publisher= {Amer. Math. Soc.},
   place= {Providence, RI},
   date= {1997},
}
  
\bib{OvSi2}{article}{
   author={Ovchinnikov, Yu.N.},
   author={Sigal, I.M.},
   title={The energy of Ginzburg-Landau vortices},
   journal={European J. Appl. Math.},
   volume={13},
   date={2002},
   number={2},
   pages={153--178},
}


\bib{OvSi3}{article}{
   author={Ovchinnikov, Yu.N.},
   author={Sigal, I.M.},
   title={Symmetry-breaking solutions of the Ginzburg-Landau equation},
   journal={J. Exp. Theor. Phys.},
   volume={99},
   date={2004},
   number={5},
   pages={1090--1107},
}
 
\bib{PaRi}{book}{
   author={Pacard, F.},
   author={Rivi\`ere, T.},
   title={Linear and nonlinear aspects of vortices. The Ginzburg-Landau model},
   series={Progress in Nonlinear Differential Equations and their Applications},
   volume={39},
   publisher={Birkh\"auser},
   place={Boston},
   date={2000},
   pages= {x+342},
}
 

\bib{Str1}{article}{
   author={Struwe, M.},
   title={On the asymptotic behavior of minimizers of the Ginzburg-Landau model in
$2$ dimensions},
   journal={Differential Integral Equations},
   volume={7},
   date={1994},
   number={5-6},
   pages={1613--1624},
}
 
\bib{Str2}{article}{
   author={Struwe, M.},
   title={Erratum: ``On the asymptotic behavior of minimizers of the Ginzburg-Landau model in
$2$ dimensions"},
   journal={Differential Integral Equations},
   volume={8},
   date={1995},
   number={1},
   pages={224},
}

\end{thebibliography}
\end{document}